\documentclass[letterpaper, 10 pt, conference]{ieeeconf}  

\usepackage[left=0.75in,right=0.75in,top=0.72in,bottom=0.78in]{geometry}
\usepackage{booktabs}
\usepackage{graphicx}
\IEEEoverridecommandlockouts                              

\usepackage{amsmath}

\usepackage{amssymb}
\usepackage{comment}
\usepackage{color,soul}
\usepackage{mwe}
\usepackage{svg}
\usepackage[hidelinks]{hyperref}
\usepackage{mathtools}
\usepackage{balance}

\usepackage{cuted}

\usepackage{acronym}
\acrodef{OFO}{Online Feedback Optimization}

 \usepackage{tikz}
\usetikzlibrary{plotmarks}
\usetikzlibrary{shapes,snakes}
\usepackage{pgfplots}

\usepackage[update,prepend]{epstopdf} 
\usepackage{fontawesome} 
\usepackage{circuitikz} 
\usepackage{tikz}
\usetikzlibrary{arrows.meta, decorations.markings}
\usepackage{cprotect} 
\usepackage[super]{nth}
\usepackage{graphicx}
\usepackage{subcaption}
\usepackage{algorithm}
\usepackage{algpseudocode}
\usepackage{listings}
\usepackage{placeins}

\usepackage{psfrag}
\usepackage{pstool}

\tikzstyle{block} = [draw, rectangle, 
    minimum height=3em, minimum width=5em]
\tikzstyle{sum} = [draw, circle, node distance=1cm]
\tikzstyle{input} = [coordinate]
\tikzstyle{output} = [coordinate]
\tikzstyle{tmp} = [coordinate]
\tikzstyle{pinstyle} = [pin edge={to-,thin,black}]

\makeatletter
\let\save@mathaccent\mathaccent
\newcommand*\if@single[3]{%
  \setbox0\hbox{${\mathaccent"0362{#1}}^H$}%
  \setbox2\hbox{${\mathaccent"0362{\kern0pt#1}}^H$}%
  \ifdim\ht0=\ht2 #3\else #2\fi
  }
\newcommand*\rel@kern[1]{\kern#1\dimexpr\macc@kerna}
\newcommand*\wideaccent[2]{\@ifnextchar^{{\wide@accent{#1}{#2}{0}}}{\wide@accent{#1}{#2}{1}}}
\newcommand*\wide@accent[3]{\if@single{#2}{\wide@accent@{#1}{#2}{#3}{1}}{\wide@accent@{#1}{#2}{#3}{2}}}
\newcommand*\wide@accent@[4]{%
  \begingroup
  \def\mathaccent##1##2{%
    \let\mathaccent\save@mathaccent
    \if#42 \let\macc@nucleus\first@char \fi
    \setbox\z@\hbox{$\macc@style{\macc@nucleus}_{}$}%
    \setbox\tw@\hbox{$\macc@style{\macc@nucleus}{}_{}$}%
    \dimen@\wd\tw@
    \advance\dimen@-\wd\z@
    \divide\dimen@ 3
    \@tempdima\wd\tw@
    \advance\@tempdima-\scriptspace
    \divide\@tempdima 10
    \advance\dimen@-\@tempdima
    \ifdim\dimen@>\z@ \dimen@0pt\fi
    \rel@kern{0.6}\kern-\dimen@
    \if#41
      #1{\rel@kern{-0.6}\kern\dimen@\macc@nucleus\rel@kern{0.4}\kern\dimen@}%
      \advance\dimen@0.4\dimexpr\macc@kerna
      \let\final@kern#3%
      \ifdim\dimen@<\z@ \let\final@kern1\fi
      \if\final@kern1 \kern-\dimen@\fi
    \else
      #1{\rel@kern{-0.6}\kern\dimen@#2}%
    \fi
  }%
  \macc@depth\@ne
  \let\math@bgroup\@empty \let\math@egroup\macc@set@skewchar
  \mathsurround\z@ \frozen@everymath{\mathgroup\macc@group\relax}%
  \macc@set@skewchar\relax
  \let\mathaccentV\macc@nested@a
  \if#41
    \macc@nested@a\relax111{#2}%
  \else
    \def\gobble@till@marker##1\endmarker{}%
    \futurelet\first@char\gobble@till@marker#2\endmarker
    \ifcat\noexpand\first@char A\else
      \def\first@char{}%
    \fi
    \macc@nested@a\relax111{\first@char}%
  \fi
  \endgroup
}
\makeatother

\newcommand\doubleoverline[1]{\overline{\overline{#1}}}

\newcommand\widebarbar{\wideaccent\doubleoverline}

\newcommand{\T}{^\mathsf{T}} 
\renewcommand{\d}{\mathrm{d}} 
\newcommand{\R}{\mathds{R}} 

\DeclareMathOperator*{\argmin}{arg\,min}

\usepackage{dsfont}

\title{\LARGE \bf Sensitivity of Online Feedback Optimization to time-varying parameters} 

\author{Marta Zagorowska$^{1,2}$, Lars Imsland$^{1}$
\thanks{*Research supported by  Marie Curie Horizon Postdoctoral Fellowship project RELIC (grant no 101063948)}
\thanks{$^{1}$ Department of Engineering Cybernetics, Norwegian University of Science and Technology, email: 
        {\tt\small lars.imsland@ntnu.no}}%
\thanks{$^{2}$ Currently with Delft Center for Systems and Control, TU Delft, email: 
        {\tt\small m.a.zagorowska@tudelft.nl}}%
}

\begin{document}

\maketitle
\thispagestyle{empty}
\pagestyle{empty}

\begin{abstract}
Online Feedback Optimization uses optimization algorithms as dynamic systems to design optimal control inputs. The results obtained from Online Feedback Optimization depend on the setup of the chosen optimization algorithm. In this work we analyse the sensitivity of Online Feedback Optimization to the parameters of projected gradient descent as the algorithm of choice. We derive closed-form expressions for sensitivities of the objective function with respect to the parameters of the projected gradient and to time-varying model mismatch. The formulas are then used for analysis of model mismatch in a gas lift optimization problem. The results of the case study indicate that the sensitivity of Online Feedback Optimization to the model mismatch depends on how long the controller has been running, with decreasing sensitivity to mismatch in individual timesteps for long operation times. 
\end{abstract}

\section{Introduction}
\label{sect:intro}
The main idea of Online Feedback Optimization relies on using optimization algorithms as dynamic systems to find the optimal control inputs. The performance of OFO controllers depends on how the parameters of the chosen optimization algorithm are set up. However, the parameters of the optimization algorithm are typically separate from the controlled system, making analysis of their impact challenging. To facilitate understanding of the impact of parameters on the controlled system, we provide an analysis of the sensitivity of OFO to parameters of projected gradient descent used as the algorithm of choice, as well as to time-varying model mismatch.

Sensitivity of classic optimization-based controllers, such as Model Predictive Control (MPC), allows shaping the performance of the controller. Sensitivity to MPC parameters has been provided in \cite{Differentiable_Amos2018} and the results were used to find optimal controller parameters. The sensitivities of MPC with respect to uncertain model parameters were used to find a robust controller in \cite{Sensitivity_Thombre2021}. Extending the robustness and tuning capabilities to Online Feedback Optimization by using sensitivity analysis is the main motivation for this work. 

The sensitivity of OFO from the perspective of robustness to model mismatch was demonstrated experimentally \cite{Non_Haeberle2020,Experimental_Ortmann2020}. A theoretical analysis of robustness for OFO with relaxed output constraints was provided in \cite{Towards_Colombino2019}. A more general approach to sensitivity of optimization-based controllers was adopted in \cite{Regularity_Mestres2023} which indicated conditions for OFO to be continuously differentiable, without providing explicit formulas for the sensitivity. In this paper, we use the fact that OFO with projected gradient descent relies on solving a series of convex quadratic optimization problems to derive closed form expressions for sensitivity.

 The contributions of this paper are:
\begin{itemize}
    \item Derivation of closed-form expressions for sensitivity of Online Feedback Optimization to its parameters;
    \item Analysis of sensitivity of OFO with respect to model mismatch for a gas lift system.
\end{itemize}

The paper is structured as follows. Online Feedback Optimization is introduced in Section \ref{sec:FO}. Section \ref{sec:Sensitivity} presents the closed-formed expressions for sensitivity of OFO, which are then numerically validated in Section \ref{sec:Numerical} in a gas lift optimization problem. Section \ref{sec:Discussion} discusses the impact of the sensitivities while Section \ref{sect:concl} indicates future directions.

\subsection{Notation}
This paper follows the notational convention for vector and matrix derivatives from \cite{Matrix_Magnus2019}. In particular, we have $\frac{\partial \mathbf{y}}{\partial x}\in\R^{n_y}$ if $ \mathbf{y}\in\R^{n_y}$ and $\frac{\partial y}{\partial \mathbf{x}}\in\R^{1\times n_x}$ if $ \mathbf{x}\in\R^{n_x}$.

\section{Online Feedback Optimization}
\label{sec:FO}
The main idea of \ac{OFO} is to treat optimization algorithms as dynamic systems. The dynamic system representing the optimization algorithm is then connected in a closed loop with the controlled system with inputs $u$ and outputs $y$. \ac{OFO} iteratively updates the input of a system $u$ to make the system converge to a local optimum of an optimization problem. A review of OFO controllers was done in \cite{Optimization_Hauswirth2021}. 

\subsection{Problem statement}
The optimization problem in \ac{OFO} is formulated as:
\begin{subequations} \label{eqn:ProblemStatement}
\begin{align}
\min_{u,y}& \quad \Phi(u,y)
    \label{eqn:CostFcn}\\
\text{subject to }    & y=h(u)\label{eq:mapping}\\
    &u\in\mathcal{U}, y\in\mathcal{Y}
\end{align}
\end{subequations}
where $\Phi:\R^{n_u}\times\R^{n_y}\rightarrow\R$ is a continuously differentiable cost function, $h:\R^{n_u}\rightarrow\R^{n_y}$ is a continuously differentiable nonlinear input-output mapping, and $\mathcal{U}$ and $\mathcal{Y}$ describe the constraints on the inputs and outputs:
\begin{align}
\mathcal{U}=&{}\lbrace u\in\R^{n_u}:Au\leq b \rbrace 
\label{eq:UBounds}\\
    \mathcal{Y}=&{}\lbrace y\in\R^{n_y}:Cy\leq d \rbrace
    \label{eq:YBounds}
\end{align}
where $A\in \R^{n_{c_1}\times n_u}$, $b\in\R^{n_{c_1}}$, $C\in\R^{n_{c_2}\times n_y}$, and $d\in\R^{n_{c_2}}$ are constant matrices \cite{Non_Haeberle2020}, and $c_1$ and $c_2$ denote the number of input and output constraints, respectively.

The \ac{OFO} controller used in this paper was proposed in \cite{Non_Haeberle2020} and uses projected gradient descent as the optimization algorithm, with a constant step size \(\alpha>0\):
\begin{align}\label{eqn:Verena_feedback}
    u^{k+1} = u^k + \alpha\widehat{\sigma}_\alpha (u^k,y^k) \quad \text{with } y^k = h(u^k),
\end{align} 
where $y^k = h(u^k)$ is the measured system output at time $k$, and \(\widehat{\sigma}_\alpha (u^k,y^k)\) is the minimizer of the constrained optimization problem (superscript $^k$ dropped for space reasons):
\begin{subequations} \label{eqn:Verena_opt}
\begin{align}
    \widehat{\sigma}_\alpha (u,y) = &\argmin_{w\in\mathbb{R}^{n_u}}\left|\left| w + G^{-1}(u)H^\top(u)\nabla\Phi^\top(u,y)\right|\right|_{G(u)}^2
    \label{eqn:Verena_sigma}\\
    &\text{subject to}\quad A\left(u+\alpha w\right)\leq b \label{eq:UProj}\\
    &\qquad\qquad\quad C \left(y+\alpha\nabla h(u) w\right)\leq d \label{eq:YProj}
\end{align}
\end{subequations}
where \(H(u)^\top = \left[\mathbb{I}_{n_u} \ \nabla h(u)^\top\right]\), $w\in\R^{n_u}$, and $G:\mathcal{U}\rightarrow\mathbb{S}_+^{n_u}$ is a continuous metric on $\mathcal{U}$. The set $\mathbb{S}_+^{n_u}$ denotes a symmetric positive definite matrix of size $n_u\times n_u$. The matrix $\mathbb{I}_{n_u}$ is an identity matrix of size $n_u\times n_u$. 

\subsection{OFO as a function of parameters}
\label{sec:OFOAsMismatchFcn}
We want to analyse how the objective $\Phi$ changes if the parameters of the system change. For simplicity, we introduce $\mathbf{p}=[p_s]_{s=0,\ldots,k}$, $p_s\in\R^{n_p}$, as the parameter whose changes affect the objective at time $k$. We have
\begin{equation}
    \Phi^k:=\Phi(u^k(\mathbf{p}),y^k(\mathbf{p}),\mathbf{p})
\end{equation}
The parameters $\mathbf{p}$ may include the elements of matrices $A$, $b$, $C$, $d$, as well as parameters of the derivatives $\nabla h$ and $\nabla\Phi$. Usually, $G$ in \eqref{eqn:Verena_opt} is chosen as a constant matrix and will be treated as a parameter together with $\alpha$.

To analyse the impact of parameters, we rewrite the projected gradient descent from \eqref{eqn:Verena_opt} as a quadratic optimization problem with linear constraints (QP). At every iteration $k$, the constraints \eqref{eq:UProj} and \eqref{eq:YProj} are linear with respect to $w$:
\begin{equation}
\label{eq:MismatchCstrRewrite}
        \underbrace{\begin{bmatrix}
        A\alpha\\ C\alpha \nabla h(u^k)
    \end{bmatrix}}_{:=\widebarbar{A}}w\leq \underbrace{\begin{bmatrix}
        b-Au^k\\ d-Cy^k
    \end{bmatrix}}_{:=\widebarbar{b}}
\end{equation}
with $\widebarbar{A}\in\R^{\widebarbar{n}\times n_u}$ and $\widebarbar{b}\in\R^{\widebarbar{n}}$, where $\widebarbar{n}=n_{c_1}+n_{c_2}$ denotes the total number of constraints in the optimization problem. For instance, if there are $n_u$ inputs and $n_y$ outputs, all constrained by bound constraints, then $\widebarbar{n}=2(n_u+n_y)$. Then we rewrite the objective function \eqref{eqn:Verena_sigma} to get a quadratic formulation \cite{Numerical_Nocedal1999}:
\begin{equation}
\label{eq:MismatchRewrittenObj}
    \begin{aligned}
        \left|\left| w + G^{-1}H^\top\nabla\Phi^\top\right|\right|_{G}^2 &{}
        =\frac{1}{2}w^\top \widebarbar{G} w+ w^\top\widebarbar{c}+\widebarbar{M}
        \end{aligned}
\end{equation}
where $\widebarbar{G}:=2G\in\R^{n_u\times n_u}$ preserves the positive definite and symmetric properties of $G$, $\widebarbar{c}:=2H^\top\nabla\Phi^\top \in\R^{n_u}$. We note that $\widebarbar{M}:=(G^{-1}H^\top\nabla\Phi^\top)^\top H^\top\nabla\Phi^\top$ is constant and independent of $w$ so it can be omitted from the problem formulation and \eqref{eqn:Verena_opt} becomes
\begin{subequations}
\label{eq:MismatchRewritten}
\begin{align}
    &\min_{w\in\mathbb{R}^{n_u}} \frac{1}{2}w^\top \widebarbar{G} w+ w^\top\widebarbar{c}
    \\
    &\text{subject to}\quad \widebarbar{A}w\leq \widebarbar{b}\label{eq:MismatchRewrittenCstr}
\end{align}
\end{subequations}
where the arguments $u$, $y$ are omitted for space reasons.

\subsection{Lagrangian}
The Lagrangian of \eqref{eq:MismatchRewritten} is:
\begin{equation}
    \mathcal{L}(w)=\frac{1}{2}w^\top \widebarbar{G} w+ w^\top\widebarbar{c}+\sum\limits_{j=1}^{\widebarbar{n}}\lambda_j(\widebarbar{a}_jw-\widebarbar{b}_j)
\end{equation}
where $\widebarbar{a}_j$ is the $j$-th row of the matrix $\widebarbar{A}$, $\widebarbar{b}_j$ is the $j$-th element of $\widebarbar{b}$, and $\lambda_j\geq 0$, $j=1,\ldots,\widebarbar{n}$ are dual variables. The matrix $\widebarbar{G}$ is positive definite, thus the problem \eqref{eq:MismatchRewritten} is convex and the KKT conditions give optimality:
\begin{subequations}
\label{eq:MismatchKKTProjection}
    \begin{align}
        \nabla_w\mathcal{L}(w)=&{}\widebarbar{G}w+\widebarbar{c}+\sum\limits_{j=1}^{\widebarbar{n}}\lambda_j\widebarbar{a}^\top_j=0\\
        \widebarbar{a}_iw-\widebarbar{b}_i\leq&{} 0,\forall i=1,\ldots,\widebarbar{n}\\
        \lambda_i(\widebarbar{a}_iw-\widebarbar{b}_i)=&{}0,\forall i=1,\ldots,\widebarbar{n}\\
        \lambda_i\geq &{} 0,\forall i=1,\ldots,\widebarbar{n}
    \end{align}
\end{subequations}
For the remainder of the paper we will assume that the QP in \eqref{eq:MismatchRewritten} is non-degenerate in the sense introduced by \cite{Sensitivity_Boot1963} and thus \eqref{eq:MismatchKKTProjection} has a solution. Following \cite{Amos_2017}, we will use \eqref{eq:MismatchKKTProjection} to find the sensitivity of OFO.

\section{Sensitivity of Online Feedback Optimization}
\label{sec:Sensitivity}
\subsection{Sensitivity to problem parameters}
\label{sec:ProblemParameters}
The matrices $\widebarbar{A}$, $\widebarbar{b}$, $\widebarbar{c}$, and $\widebarbar{G}$ change with iterations $s$ because they depend on current measurements $u^s$, $y^s$, or on the parameters $p_s$. At time $k$ we want to find:
\begin{equation}
\label{eq:PhiParamDeriv}
\begin{aligned}
     \frac{\partial\Phi^k}{\partial p_s} =&{} \frac{\partial\Phi^k}{\partial p_s} + \left(\frac{\partial\Phi^k}{\partial u^{k}}+\frac{\partial\Phi^k}{\partial y^{k}}\cdot\frac{\partial y^k}{\partial u^{k}}\right)\cdot\frac{\partial u^k}{\partial p_{s}} \in\R^{1\times n_p}
    \end{aligned}
\end{equation}
for $s=0,\ldots,k$. From \eqref{eqn:Verena_feedback}, we see that $u^k$ and $y^k$ depend only on $p_s$ for $s=0,\ldots,k-1$, and hence:
\begin{equation}
    \left(\frac{\partial\Phi^k}{\partial u^{k}}+\frac{\partial\Phi^k}{\partial y^{k}}\cdot\frac{\partial y^k}{\partial u^{k}}\right)\cdot\frac{\partial u^k}{\partial p_{k}}=0
\end{equation}
in \eqref{eq:PhiParamDeriv}. As a result, $\frac{\partial\Phi^k}{\partial p_k}$ corresponds to computing instantaneous sensitivity of $\Phi^k$ with respect to parameter $p^k$, which is unaffected by OFO and will be omitted from further analysis.

We also see that $\frac{\partial\Phi^k}{\partial u^k}+\frac{\partial\Phi^k}{\partial y^k}\cdot \frac{\partial y^k}{\partial u^k}$ in \eqref{eq:PhiParamDeriv} is independent of $s$ and thus we can write in matrix form:
\begin{equation}
    \begin{bmatrix}
        \frac{\partial\Phi^k}{\partial p_{k-1}}\\
        \frac{\partial\Phi^k}{\partial p_{k-2}}\\
        \vdots\\
        \frac{\partial\Phi^k}{\partial p_{0}}\\
    \end{bmatrix}=\mathbb{I}_{k}\otimes \left(\frac{\partial\Phi^k}{\partial u^k}+\frac{\partial\Phi^k}{\partial y^k}\cdot \frac{\partial y^k}{\partial u^k}\right)\cdot     \begin{bmatrix}
        \frac{\partial u^k}{\partial p_{k-1}}\\
        \frac{\partial u^k}{\partial p_{k-2}}\\
        \vdots\\
        \frac{\partial u^k}{\partial p_{0}}\\
    \end{bmatrix}
    \label{eq:dpartialPhi}
\end{equation}
where $\otimes$ denotes the Kronecker product and $\mathbb{I}_{k}$ is an identity matrix $k\times k$.

From \eqref{eq:dpartialPhi} we obtain the total derivative of $\Phi$ at time $k$:
\begin{equation}
\begin{aligned}  
   \text{d} \Phi^k =&{}\sum\limits_{i=0}^{k-1} \frac{\partial\Phi^k}{\partial p_{i}}\cdot\text{d} p_i
    \end{aligned}
\end{equation}
where $\Delta p_{i}$ are individual increments at time $i$.

\subsection{Derivative of $u$}
We see in \eqref{eq:dpartialPhi} that we need $\frac{\partial u^k}{\partial p_s}\in\R^{n_u\times n_p}$, $s=0,\ldots,k-1$. 
We have from \eqref{eqn:Verena_feedback}:
\begin{equation}
    u^k=u^{k-1}(p_0,\ldots,p_{k-2})+\alpha w^{k-1}_*(p_0,\ldots,p_{k-1})
\end{equation}
We assume for now that $\alpha$ and $u_0$ are independent of $\mathbf{p}$, thus:
\begin{subequations}
\begin{align}
\label{eq:RewrittenBigDerivativeU}
    \frac{\partial u^k}{\partial p_s}=&{}\frac{\partial u^{k-1}}{\partial p_s}+\frac{\partial \left( \alpha w^{k-1}_*\right)}{\partial p_s}\\
    =&{} \alpha\sum\limits_{i=s}^{k-1}\frac{\partial w^{i}_*}{\partial p_s}
\end{align}
\end{subequations}
for $s=0,\ldots,k-1$. In matrix form, we obtain:
\begin{equation}
    \begin{bmatrix}
        \frac{\partial u^k}{\partial p_{k-1}}\\
        \frac{\partial u^k}{\partial p_{k-2}}\\
        \vdots\\
        \frac{\partial u^k}{\partial p_{0}}\\
    \end{bmatrix}=\alpha \begin{bmatrix}
                \frac{\partial w^{k-1}}{\partial p_{k-1}}&0&\\
        \frac{\partial w^{k-1}}{\partial p_{k-2}}&\frac{\partial w^{k-2}}{\partial p_{k-2}}&0\\
        \vdots&&\vdots\\
        \frac{\partial w^{k-1}}{\partial p_{0}}&\frac{\partial w^{k-2}}{\partial p_{0}}&\frac{\partial w^0}{\partial p_{0}}\\
    \end{bmatrix}\cdot \begin{bmatrix}
        \mathbb{I}_{n_p} \\\vdots\\  \mathbb{I}_{n_p}
    \end{bmatrix}
    \label{eq:finalderivu}
\end{equation}

Thus, we obtain a total derivative $\text{d}u^k$ as:
\begin{equation}
\begin{aligned}
    \text{d}u^k=&{} \alpha\sum\limits_{s=0}^{k-1}\sum\limits_{i=s}^{k-1}  \frac{\partial w^{i}}{\partial p_{s}} \text{d} p_{s}
    \end{aligned}
    \label{eq:totalu}
\end{equation}
The result from \eqref{eq:totalu} corresponds to the formulas provided by \cite{Backpropagation_Werbos1990} for backpropagation in time.

\subsection{Derivative of $w$}
To obtain the derivatives of $u$ from \eqref{eq:totalu}, we need to calculate $\frac{\partial w^{i}}{\partial p_{s}}\in\R^{n_u\times n_p}$. Taking into account that at time $k$ the matrices $\widebarbar{A}$, $\widebarbar{b}$, $\widebarbar{c}$, $\widebarbar{G}$ depend on the inputs and the outputs up to $s=k-1$, and the parameters up to $s=k$, we  can compute $\frac{\partial w_*^k}{\partial p_{s}}$ using the chain rule:
\begin{equation}
\label{eq:BigWDerivative}
\begin{aligned}
    \frac{\partial w^k_*}{\partial p_{s}}=&{}\underbrace{K^s}_{=0 \text{ for } s=0,\ldots,k-1}+
    \\&{}\underbrace{\mathrm{d}w_{\widebarbar{b}}\mathrm{d}_{p_{s}}\widebarbar{b}+\mathrm{d}w_{\widebarbar{A}}\mathrm{d}_{p_{s}}\widebarbar{A}+\mathrm{d}w_{\widebarbar{c}}\mathrm{d}_{p_{s}}\widebarbar{c}+\mathrm{d}w_{\widebarbar{G}}\mathrm{d}_{p_{s}}\widebarbar{G}}_{=0 \text{ for } s=k}
    \end{aligned}
\end{equation}
where $K^k = \mathrm{d}w_{\widebarbar{b}}+\mathrm{d}w_{\widebarbar{A}}+\mathrm{d}w_{\widebarbar{c}}+\mathrm{d}w_{\widebarbar{G}}$. From \cite{Amos_2017}, the terms $\mathrm{d}w_{\widebarbar{b}}$, $\mathrm{d}w_{\widebarbar{A}}$, $\mathrm{d}w_{\widebarbar{c}}$, and $w_{\widebarbar{G}}$ are obtained by differentiating \eqref{eq:MismatchKKTProjection}:
\begin{equation}
\begin{aligned}
    \begin{bmatrix}
        \widebarbar{G}& \widebarbar{A}^\top\\
        D(\lambda^*)\widebarbar{A}&D(\widebarbar{A}w^*-\widebarbar{b})
    \end{bmatrix}\cdot \begin{bmatrix}
        \d w_X\\ \d \lambda_X
    \end{bmatrix}=\\-
    \begin{bmatrix}
        \d  \widebarbar{G}w^*+\d \widebarbar{c}+\d \widebarbar{A}^\top\lambda^*\\
         D(\lambda^*)\d\widebarbar{A}w^*-D(\lambda^*)\d \widebarbar{b}
    \end{bmatrix}
    \end{aligned}
\end{equation}
where $D(x)$ creates a diagonal matrix with entries $x$, $X\in\lbrace \widebarbar{b},\widebarbar{A},\widebarbar{c},\widebarbar{G}\rbrace$. For instance if we want to find $\d w_{\widebarbar{b}}\in\R^{n_u\times \widebarbar{n}}$, we set the required element of $\mathrm{d}\widebarbar{b}$ to unity, $\mathrm{d}\widebarbar{A}=\mathrm{d}\widebarbar{c}=\mathrm{d}\widebarbar{G}=0$ and solve the resulting equations for $\d w_{\widebarbar{b}}$. 

We now need to compute derivatives of $\widebarbar{A}$, $\widebarbar{b}$, $\widebarbar{c}$, $\widebarbar{G}$ with respect to the parameters $p_s$. For instance,  we have:
\begin{equation}
\begin{aligned}
    \mathrm{d}_{p_{s}}&{}\widebarbar{b}(u^k,y^k) = \underbrace{\frac{\partial\widebarbar{b}(u^k,y^k)}{\partial p_{s}}}_{=0 \text{ for } s=0,\ldots,k-1}+\\
    &{}\underbrace{\left(\frac{\partial\widebarbar{b}(u^k,y^k)}{\partial u_{k}}+\frac{\partial\widebarbar{b}(u^k,y^k)}{\partial y_{k}}\cdot\frac{\partial y^k}{\partial u_{k}}\right)\cdot\frac{\partial u^k}{\partial p_{s}}}_{=0 \text{ for } s=k}\in\R^{\widebarbar{n}\times n_p}
    \end{aligned}
\end{equation}
The formulas for $\mathrm{d}_{p_{s}}\widebarbar{A}\in\R^{\widebarbar{n}\times n_un_p}$, $\mathrm{d}_{p_{s}}\widebarbar{c}\in\R^{n_u\times n_p}$, $\mathrm{d}_{p_{s}}\widebarbar{G}\in\R^{n_u\times n_un_p}$ are obtained analogously.

We can now write:
\begin{equation}
\begin{aligned}
    \begin{bmatrix}
    \frac{\partial w^k_*}{\partial p_k} \\
    \frac{\partial w^k_*}{\partial p_{k-1}}\\
    \vdots\\
    \frac{\partial w^k_*}{\partial p_0}
    \end{bmatrix}
     =&{} \begin{bmatrix}
         K^k& 0\\
         0& \mathbb{I}_{k}\otimes K
     \end{bmatrix}  \cdot \begin{bmatrix}
     \mathbb{I}_{n_u}\\
        \frac{\partial u^k}{\partial p_{k-1}}\\
        \frac{\partial u^k}{\partial p_{k-2}}\\
        \vdots\\
        \frac{\partial u^k}{\partial p_{0}}\\
    \end{bmatrix}
     \end{aligned}
     \label{eq:FinalDerivw}
\end{equation}
where $K= K_{\widebarbar{b}}+K_{\widebarbar{A}}+K_{\widebarbar{c}}+K_{\widebarbar{G}}$ and $K_X=\mathrm{d}w_{X}\left(\frac{\partial X}{\partial u_{k}}+\frac{\partial X}{\partial y_{k}}\cdot\frac{\partial y^k}{\partial u_{k}}\right)$, $X\in\lbrace \widebarbar{b},\widebarbar{A},\widebarbar{c},\widebarbar{G}\rbrace$.

Inserting \eqref{eq:FinalDerivw} into \eqref{eq:finalderivu} and subsequently into \eqref{eq:dpartialPhi}, we obtain closed-form expressions for sensitivity of the objective function $\Phi$ at time $k$ to the parameter $\mathbf{p}$.

\subsection{Sensitivity to OFO parameters}
The sensitivities derived in Section \ref{sec:ProblemParameters} allow analysing how the parameters of the original problem \eqref{eqn:ProblemStatement} affect the control input $u^k$ at time $k$ and the corresponding objective function $\Phi^k$. In a similar way to the derivations in Section \ref{sec:ProblemParameters}, we will now present sensitivities with respect to parameters of OFO: $\alpha$, $G$, and $u_0$.

\subsubsection{Sensitivity to $G$}
The sensitivity of the control input and the objective function to the matrix $G$ can be obtained directly from the formulas in Section \ref{sec:ProblemParameters} by taking $p_s=G$.

\subsubsection{Sensitivity to $\alpha$}
When computing the sensitivity with respect to $\alpha$, the formulas \eqref{eq:dpartialPhi} and \eqref{eq:FinalDerivw} remain as in Section \ref{sec:ProblemParameters}. To obtain the equivalent to formula \eqref{eq:finalderivu}, we notice from \eqref{eqn:Verena_feedback} that:
\begin{equation}
    u^{k}=u^0+\sum\limits_{s=0}^{k-1}\alpha_s w^s
    \label{eq:RecursiveU}
\end{equation}
where for the ease of notation we took $\alpha=\alpha_s$ as a time-varying parameter. Assuming that $u^0$ is independent of $\alpha_s$ and taking the derivatives of \eqref{eq:RecursiveU} at time $k$ with respect to $\alpha_s$, we get:
\begin{equation}
    \frac{\partial u^k}{\partial \alpha_s}=w^s+\sum\limits_{i=s}^{k-1}\alpha_s\frac{\partial w^i}{\partial \alpha_s}
\end{equation}
In turn, the equivalent to \eqref{eq:finalderivu} becomes:
\begin{equation}
    \begin{bmatrix}
        \frac{\partial u^k}{\partial \alpha_{k-1}}\\
        \frac{\partial u^k}{\partial \alpha_{k-2}}\\
        \vdots\\
        \frac{\partial u^k}{\partial \alpha_{0}}\\
    \end{bmatrix}=\begin{bmatrix}
        w^{k-1}\\w^{k-2}\\\vdots\\w^0
    \end{bmatrix}+\alpha \begin{bmatrix}
                \frac{\partial w^{k-1}}{\partial \alpha_{k-1}}&0_{k-1}&\\
        \frac{\partial w^{k-1}}{\partial \alpha_{k-2}}&\frac{\partial w^{k-2}}{\partial \alpha_{k-2}}&0_{k-1}\\
        \vdots\\
        \frac{\partial w^{k-1}}{\partial \alpha_{0}}&\frac{\partial w^{k-2}}{\partial \alpha_{0}}&\frac{\partial w^0}{\partial \alpha_{0}}\\
    \end{bmatrix}\cdot \begin{bmatrix}
        \mathbb{I}_{n_p} \\\vdots\\  \mathbb{I}_{n_p}
    \end{bmatrix}
    \label{eq:finalderivualpha}
\end{equation}

\subsection{Sensitivity to $u_0$}
In this subsection we assume that the starting point for OFO at time $k=0$ is $u^0=u_0$ and analyse the sensitivity to its value at time $k\geq 1$. From \eqref{eq:dpartialPhi}, we get that $\frac{\partial\Phi^k}{\partial u_{s}}=0$ for $s=1,\ldots,k-1$ and therefore:
\begin{equation}
        \frac{\partial\Phi^k}{\partial u_{0}}=\left(\frac{\partial\Phi^k}{\partial u_k}+\frac{\partial\Phi^k}{\partial y_k}\cdot \frac{\partial y^k}{\partial u_k}\right)\cdot  
        \frac{\partial u^k}{\partial u_{0}}\\
    \label{eq:dpartialPhiu0}
\end{equation}

To obtain $\frac{\partial u^k}{\partial u_{0}}$ we use \eqref{eq:RecursiveU} and assume that $\alpha$ is independent of $u_0$ to obtain:
\begin{equation}
        \frac{\partial u^k}{\partial u_0}=\mathbf{1}_{n_u}+\alpha\sum\limits_{i=0}^{k-1}\frac{\partial w^i}{\partial u_0}
            \label{eq:dpartialuu0}
\end{equation}
where $\mathbf{1}_{n_u}=[1]\in\R^{n_u}$. Finally, to obtain $\frac{\partial w^i}{\partial u_0}$ in \eqref{eq:dpartialuu0}, we notice that the matrices from \eqref{eq:MismatchRewritten} are independent of $u_0$ in a direct way, so the solution $w^k$ at time $k$ will only depend on $u_0$ through $u^{k-1}$ and $y^{k-1}$. Thus, from \eqref{eq:FinalDerivw} we get:

\begin{equation}
    \frac{\partial w^k_*}{\partial u_0}
     = K\cdot  \frac{\partial u^k}{\partial u_{0}}\\
     \label{eq:FinalDerivwu0}
\end{equation}

\section{Numerical results}
\label{sec:Numerical}
\subsection{Validation with finite differences}
\begin{figure}[!tbp]
     \centering
        \psfrag{fu}[][]{\tiny\textsf{$f(u)$}}
        \psfrag{Optimum}[][]{\tiny\textsf{Optimum}}
\includegraphics[width=0.3\textwidth]{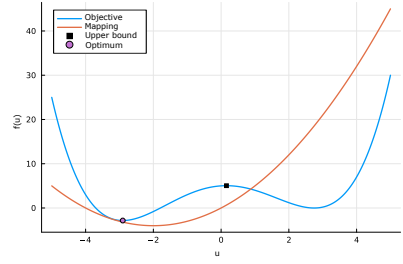}
        \caption{The objective \eqref{eq:ToyObjective} and the mapping \eqref{eq:ToyMapping} with the optimum (circle) and the upper bound of the initial condition (square)}
        \label{fig:OneDSystem}
\end{figure}

To validate the formulas obtained in Section \ref{sec:Sensitivity}, we compared the derived sensitivities with the computations based on finite differences for a one-dimensional optimization problem adapted from \cite{simulationlib} (Fig. \ref{fig:OneDSystem}):
\begin{subequations}
\begin{align}
\min_{u,y}& \quad 0.1(u^2y-4uy+5u)+5\label{eq:ToyObjective}\\
\text{subject to }    & y=u^2+4u\label{eq:ToyMapping}\\
    &u\in\mathcal{U}_i, y\in\mathcal{Y}, i=1,2
\end{align}
\end{subequations}
where $\mathcal{U}_1=\mathcal{Y}=[-5,5]$,  $\mathcal{U}_2=[-2,2]$. The validation was done for $\alpha\in [0.0001,0.025]$ with $\Delta\alpha=1e-4$ for finite differences ($G=1$, $u_0=-0.63$), $G\in[0.5,40]$ with $\Delta G=5e-2$ ($\alpha=0.01$, $u_0=-0.63$), and $u_0\in[-4,0.156]$ with $\Delta u_0=5e-3$, ($G=1$, $\alpha=0.01$), for three simulation times $T_F\in\lbrace 50, 100, 150\rbrace$.

\subsubsection{Analysis for unconstrained optimum}
\begin{figure*}[!tbp]
     \centering
     \begin{subfigure}[b]{0.3\textwidth}
         \centering
         \includegraphics[width=\textwidth]{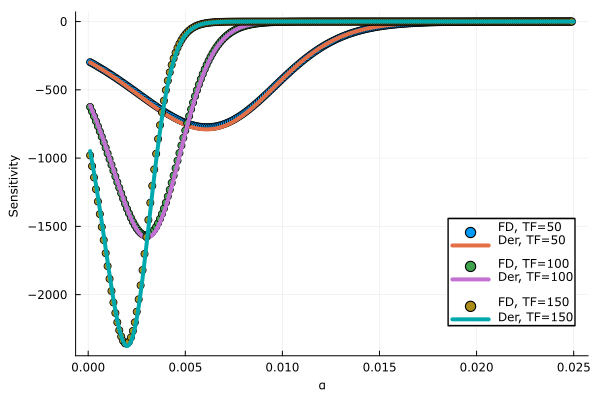}
         \caption{Sensitivity to $\alpha$, $\alpha\in[0.0001,0.025]$}
         \label{fig:dphidalpha_time}
     \end{subfigure}
     \hfill
     \begin{subfigure}[b]{0.3\textwidth}
         \centering
         \includegraphics[width=\textwidth]{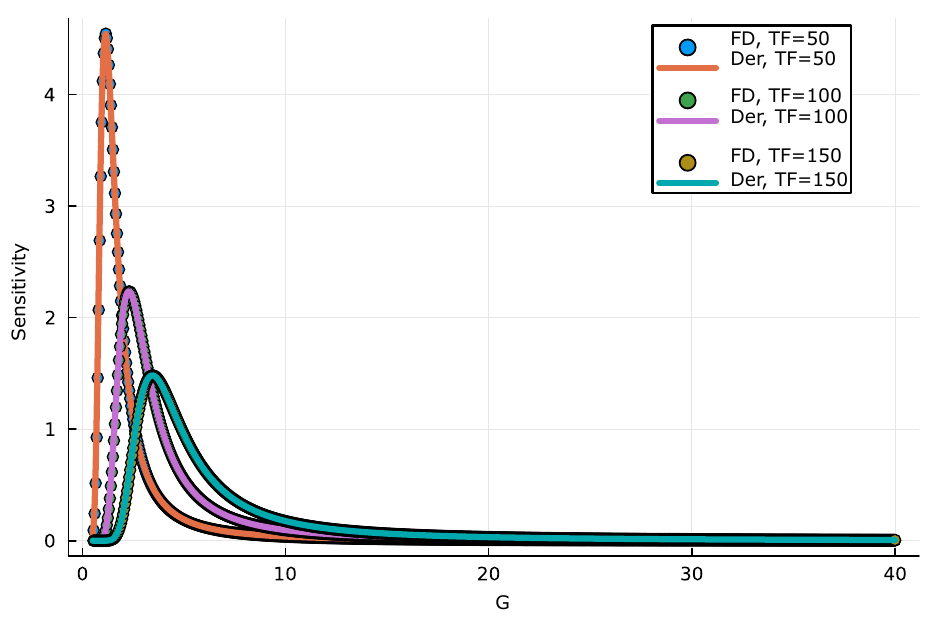}
         \caption{Sensitivity to $G$, $G\in[0.5,40]$}
         \label{fig:dphidG_alpha}
     \end{subfigure}
     \hfill
          \begin{subfigure}[b]{0.3\textwidth}
         \centering
         \includegraphics[width=\textwidth]{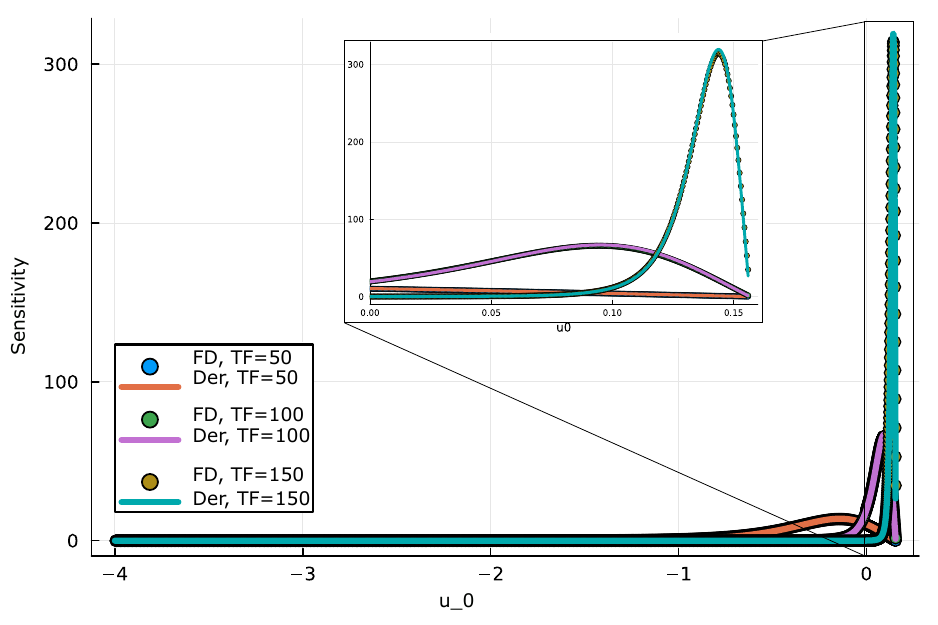}
         \caption{Sensitivity to $u_0$, $u_0\in[-4,0.156]$}
         \label{fig:dphidu0_alpha}
     \end{subfigure}

          \begin{subfigure}[b]{0.3\textwidth}
         \centering
         \includegraphics[width=\textwidth]{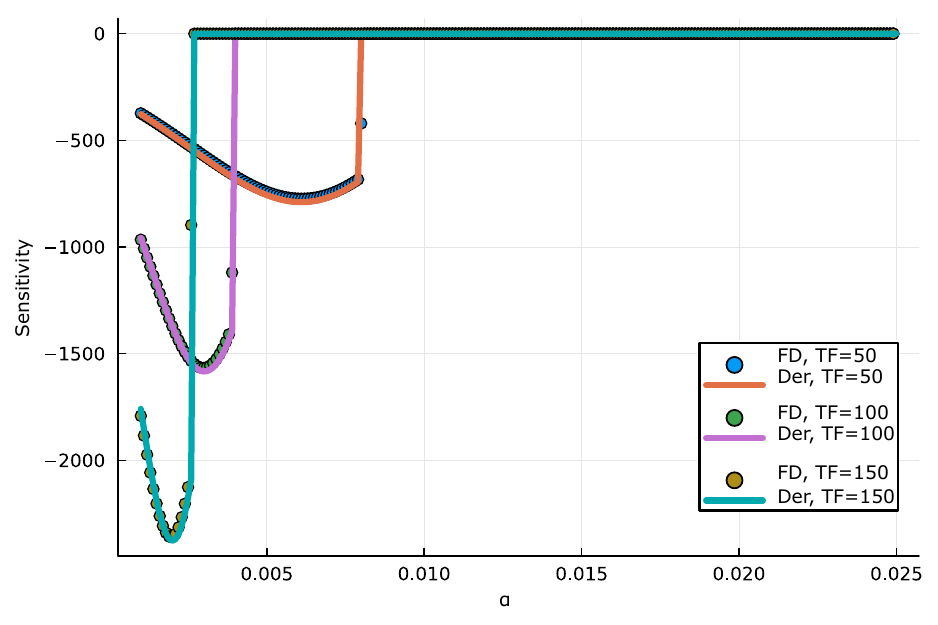}
         \caption{Sensitivity to $\alpha$, $\alpha\in[0.0001,0.025]$}
         \label{fig:dphidalpha_time_cstr}
     \end{subfigure}
     \hfill
     \begin{subfigure}[b]{0.3\textwidth}
         \centering
         \includegraphics[width=\textwidth]{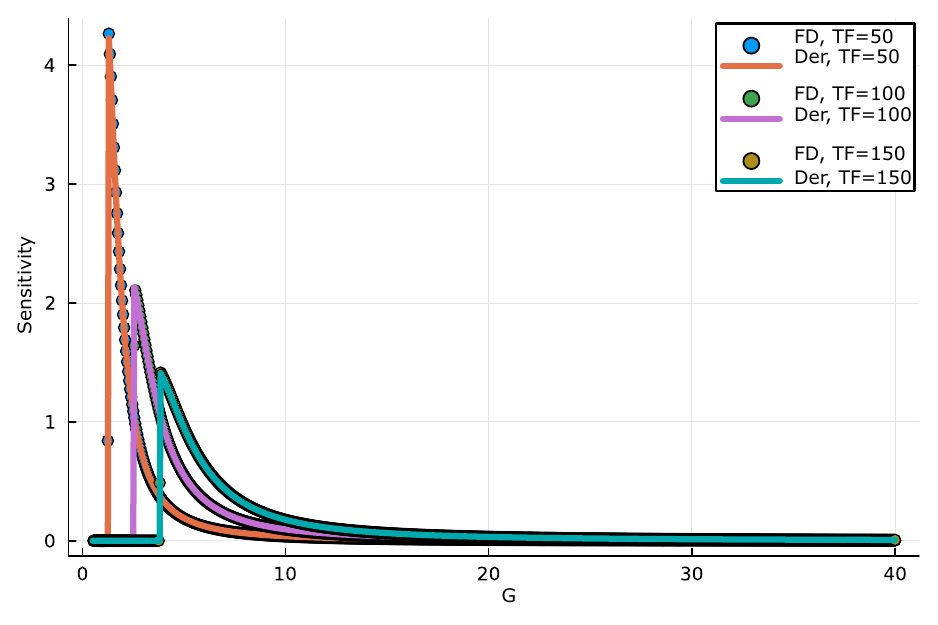}
         \caption{Sensitivity to $G$, $G\in[0.5,40]$}
         \label{fig:dphidG_alpha_cstr}
     \end{subfigure}
     \hfill
          \begin{subfigure}[b]{0.3\textwidth}
         \centering
         \includegraphics[width=\textwidth]{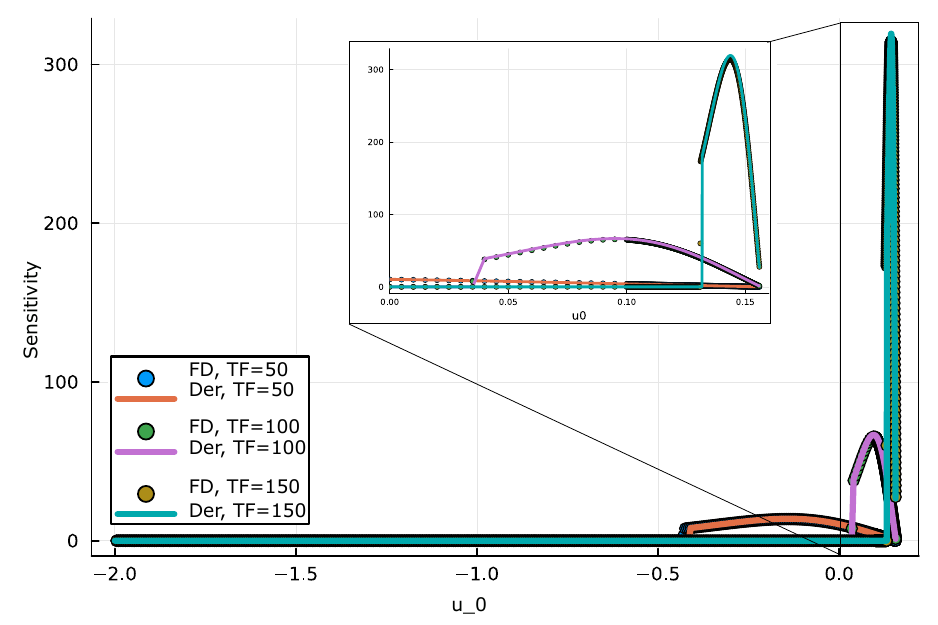}
         \caption{Sensitivity to $u_0$, $u_0\in[-2,0.156]$}
         \label{fig:dphidu0_alpha_cstr}
     \end{subfigure}
        \caption{Comparison of the formulas from Section \ref{sec:Sensitivity} (Der) to finite differences (FD) for different $T_F$}
        \label{fig:Validation}
\end{figure*}

We first validate the formulas for $i=1$ which leads to inactive constraints. Figure \ref{fig:Validation} shows a comparison of the derivatives from Section \ref{sec:Sensitivity} (solid line) and finite differences (circles). For every time $T_F$, the results obtained from the formulas in Section \ref{sec:Sensitivity} reflect the results from finite differences. This indicates the potential of using the sensitivities to locally approximate the objective as a function of the parameters, for instance for tuning purposes.

We note that the sensitivities change depending on the time horizon $T_F$. Figure \ref{fig:dphidalpha_time} indicates that the sensitivities for $\alpha\geq 0.175$ are close to zero. This is because such $\alpha$ allows reaching the optimum of \eqref{eqn:CostFcn}, driving the derivatives with respect to $u$ and $y$ to zero in \eqref{eq:PhiParamDeriv}. However, the smaller the value of $\alpha$, the more time the controller needs to reach the optimum \cite{Tuning_Zagorowska2024}. As a result, the values $\alpha\leq 0.1$ lead to high sensitivity because the controller does not reach the optimum and the derivatives of the objective at $T_F$ remain non-zero. The behaviour is further confirmed by the negative sign of the derivative in Fig. \ref{fig:dphidalpha_time}, indicating that the objective  \eqref{eqn:CostFcn} is a decreasing function of $\alpha$.

The sensitivity of the objective to the scaling parameter $G$ in Fig. \ref{fig:dphidG_alpha} indicates that the value of the objective increases with $G$. Conversely to the sensitivity to $\alpha$, $G\leq 2$ allows reaching the optimum, driving the sensitivities to zero. The relationship between $G$ and $\alpha$ has been done in \cite{Tuning_Zagorowska2024}.

Figure \ref{fig:dphidu0_alpha} indicates increased sensitivity for the initial condition close to the boundary of the region of attraction of the sought minimum. An explanation for this increased sensitivity is provided by Fig. \ref{fig:OneDSystem} and analysis of formula \eqref{eqn:Verena_opt} from the perspective of $\nabla\Phi$. The objective function in \eqref{eq:ToyObjective} (blue in Fig. \ref{fig:OneDSystem}) has a local maximum at 0.156 and thus $\nabla\Phi\approx 0$ in \eqref{eqn:Verena_opt} for $u_0$ in the neighbourhood of the maximum. As a results, the solution $w\approx 0$ yielding slow increments in $u$ from \eqref{eqn:Verena_feedback} leading to small changes in $y$ from \eqref{eq:ToyMapping}. For shorter time periods, $T_F=50$ and $T_F=100$, the slow changes in $u$ and $y$ have little impact on the objective, thus keeping $\frac{\d\Phi}{\d u_0}$ smaller. Longer time periods, $T_F=150$, allow OFO to escape the neighbourhood of the local maximum without reaching the optimum, leading to increased sensitivity of $\Phi(u_{T_F},y_{T_F})$ to $u_0$. The optimum is reached if the final time is increased, $T_F=300$, making the objective less sensitive. These results emphasise the importance of the choice of the time as well as the step size $\alpha$.

\subsubsection{Analysis for constrained optimum}
The formulas from Section \ref{sec:Sensitivity} remain valid if the constraints are active, as long as the KKT conditions \eqref{eq:MismatchKKTProjection} are nondegenerate. In particular, if the objective function \eqref{eqn:CostFcn} becomes constant when either the input or the output constraints become active, the sensitivity becomes zero (Fig. \ref{fig:dphidalpha_time_cstr}, \ref{fig:dphidG_alpha_cstr}, \ref{fig:dphidu0_alpha_cstr}).

\subsection{Usage}
Obtaining the gradients of the output mapping in \eqref{eqn:Verena_opt} may be a limitation and several works  have proposed approaches to estimate $\nabla h$ \cite{Model_He2023}. Conversely, some works focus on ensuring robustness for constant approximation of $\nabla h$ \cite{Towards_Colombino2019}. The sensitivity analysis from Section \ref{sec:ProblemParameters} allows analysing how a mismatch in $\nabla h$ over time affects the objective and the constraints. For that we look at a gas lift case study adapted from \cite{Data_Andersen2023}. 

\subsubsection{Gas lift optimization}
\label{sec:NoCoupling}
\begin{figure}[!tbp]
     \centering
         \includegraphics[width=0.3\textwidth]{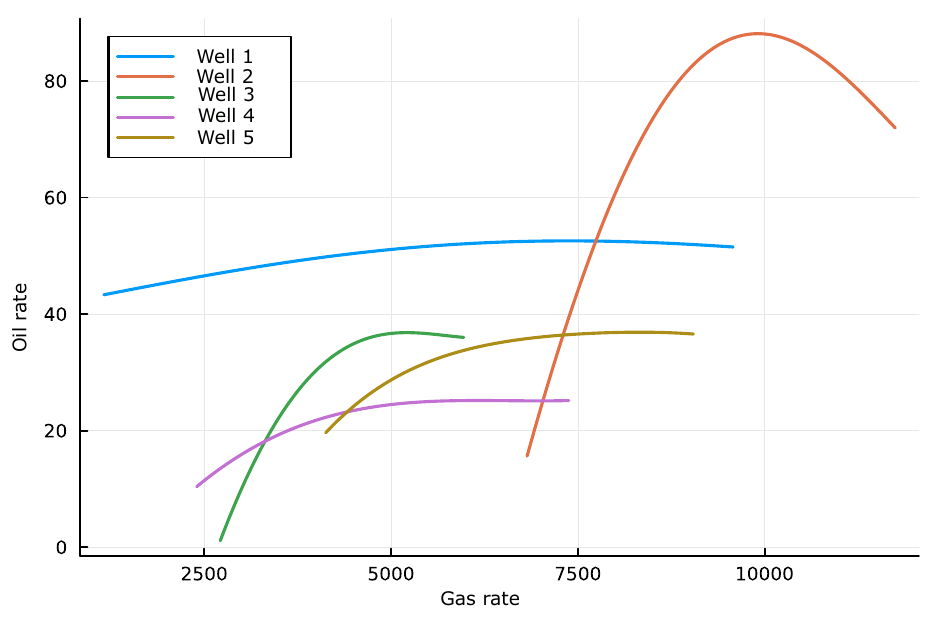}
        \caption{Characteristics of the wells adapted from \cite{Data_Andersen2023} (oil and gas rates in Sm$^3$day$^{-1}$)}
        \label{fig:Characteristics}
\end{figure}

The objective it to maximise the cumulative output of two floating oil platforms $\max_{\mathbf{u},\mathbf{y}} y_1(\mathbf{u})+y_2(\mathbf{u})
$, $\mathbf{u}=[u_i]_{i=1,\ldots,5}$,$\mathbf{y}=[y_i]_{i=1,2}$ connected to two $y_1(\mathbf{u})=f_1(\mathbf{u})+f_2(\mathbf{u})$ and three wells $y_2(\mathbf{u})= f_3(\mathbf{u})+f_4(\mathbf{u})+f_5(\mathbf{u})$ where $f_i(\mathbf{u})=\sum_{j=0}^4 s_{4-j}u_i^{4-j}$ describes the characteristics of the $i$-th well as a function of the amount of natural gas, $u_i$, injected into the well to facilitate oil extraction. The characteristics used in this work are shown in Fig. \ref{fig:Characteristics} (obtained from \cite{Data_Andersen2023} using \cite{WebPlotDigitizer_Rohatgi2018}). Both the inputs and the outputs are bounded $u_i\in[\underline{u}_i,\overline{u}_i]_{i=1,\ldots,5}$, $y_i\in[\underline{y}_i,\overline{y}_i]_{i=1,2}$.

The case study is put in the OFO framework \eqref{eqn:ProblemStatement} with $\Phi(\mathbf{u},\mathbf{y})=-y_1(\mathbf{u})-y_2(\mathbf{u})$, $h(\mathbf{u})=\begin{bmatrix}
    f_1(\mathbf{u})+f_2(\mathbf{u}) & f_3(\mathbf{u})+f_4(\mathbf{u})+f_5(\mathbf{u})
\end{bmatrix}^{\T}$, and constraint matrices:
\begin{center}
    $A=\begin{bmatrix}
          1 &  0 &  0 &  0&   0\\
 -1 &  0  & 0 &  0 &  0\\
  0  & 1 &  0 &  0 &  0\\
  0  &-1  & 0 &  0  & 0\\
  0   &0&   1 &  0 &  0\\
  0  & 0 & -1 &  0  & 0\\
  0  & 0 &  0 &  1 &  0\\
  0  & 0 &  0 & -1 &  0\\
  0 &  0 &  0 &  0&   1\\
  0  & 0 &  0   &0 & -1\\
    \end{bmatrix}$, $b=\begin{bmatrix}
        9576\\-1157\\11745\\-6819\\5972\\-2714\\7377\\-2399\\9043\\-4125
    \end{bmatrix}$, $C=\begin{bmatrix}
        1& 0\\ -1 &0\\ 0 &1\\ 0 &-1
    \end{bmatrix}$, $d=\begin{bmatrix}
        150\\0\\150\\0
    \end{bmatrix}$
\end{center}
The initial condition was set $u_0=[2500.0,7000.0,4500.0,4500.0,4500.0]$ Sm$^3$day$^{-1}$.

\begin{figure*}[!tbp]
     \centering
     \begin{subfigure}[b]{0.3\textwidth}
         \centering
         \includegraphics[width=\textwidth]{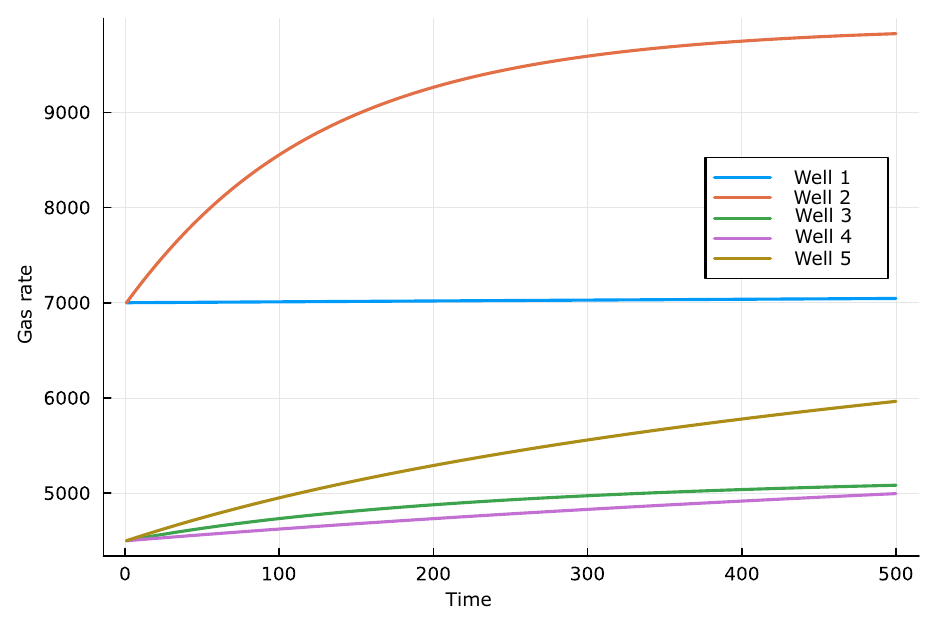}
         \caption{Inputs $u_i$ for the five wells, in Sm$^3$day$^{-1}$}
         \label{fig:Inputs1}
     \end{subfigure}
     \hfill
     \begin{subfigure}[b]{0.3\textwidth}
         \centering
         \includegraphics[width=\textwidth]{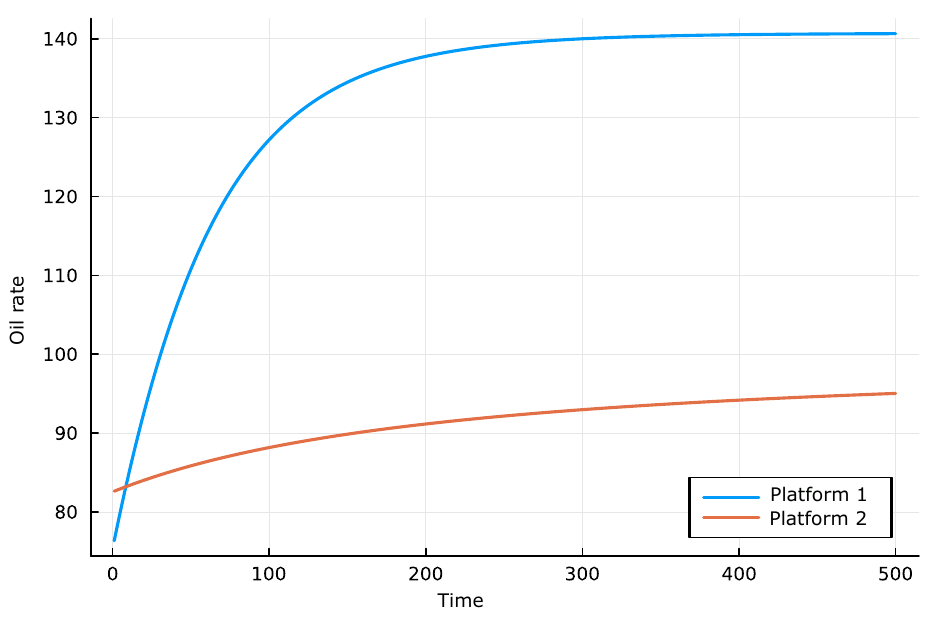}
         \caption{Outputs $y_i$ of the two platforms, in Sm$^3$day$^{-1}$}
         \label{fig:Outputs1}
     \end{subfigure}
     \hfill
          \begin{subfigure}[b]{0.3\textwidth}
         \centering
         \includegraphics[width=\textwidth]{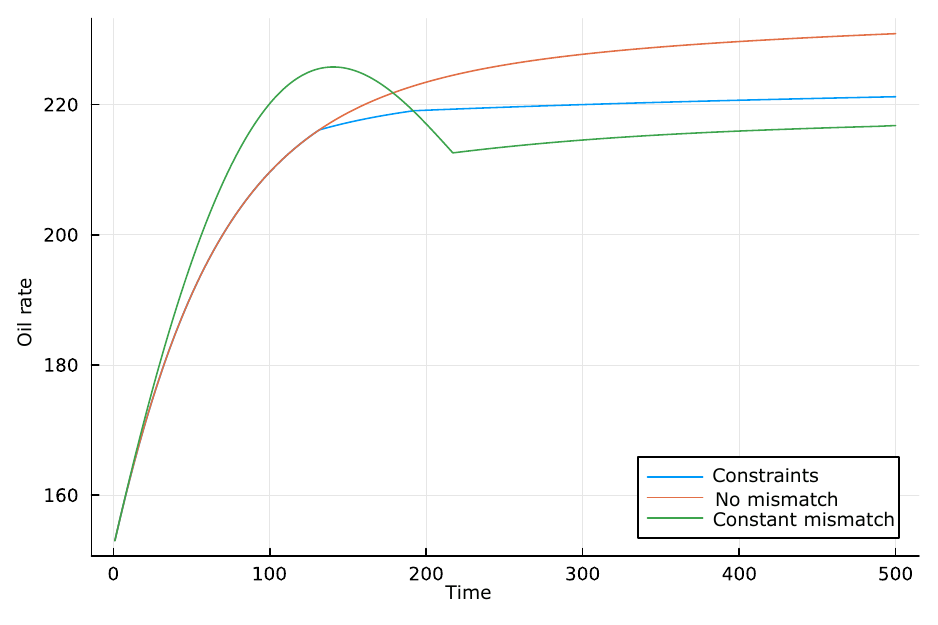}
         \caption{The objective $-\Phi$ for the three cases, in Sm$^3$day$^{-1}$}
         \label{fig:Objective}
     \end{subfigure}
     \hfill     
     \begin{subfigure}[b]{0.3\textwidth}
         \centering
         \includegraphics[width=\textwidth]{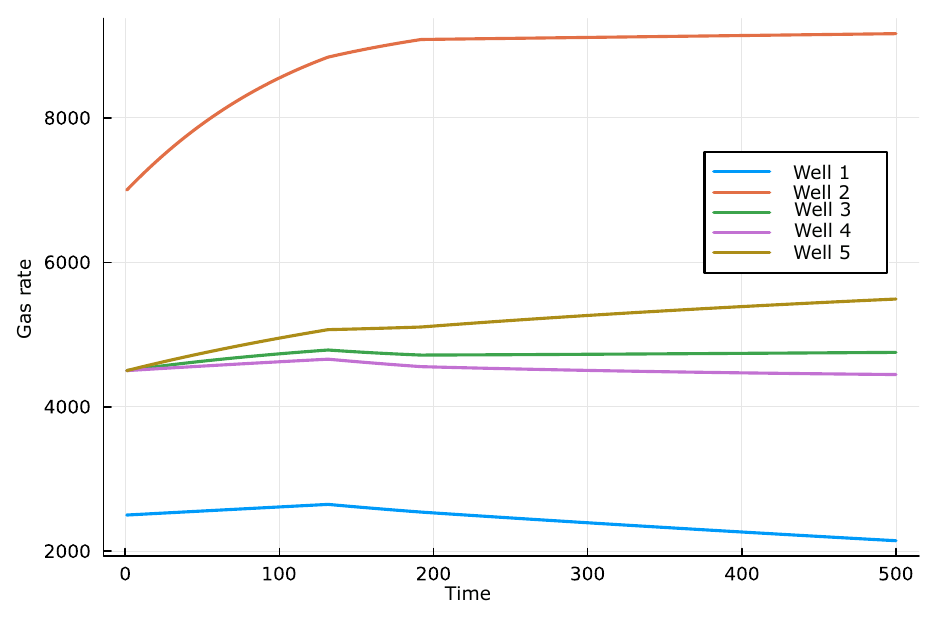}
         \caption{Inputs $u_i$ for the five wells for the constrained case, in Sm$^3$day$^{-1}$}
         \label{fig:Inputs1_cstr}
     \end{subfigure}
     \hfill
     \begin{subfigure}[b]{0.3\textwidth}
         \centering
         \includegraphics[width=\textwidth]{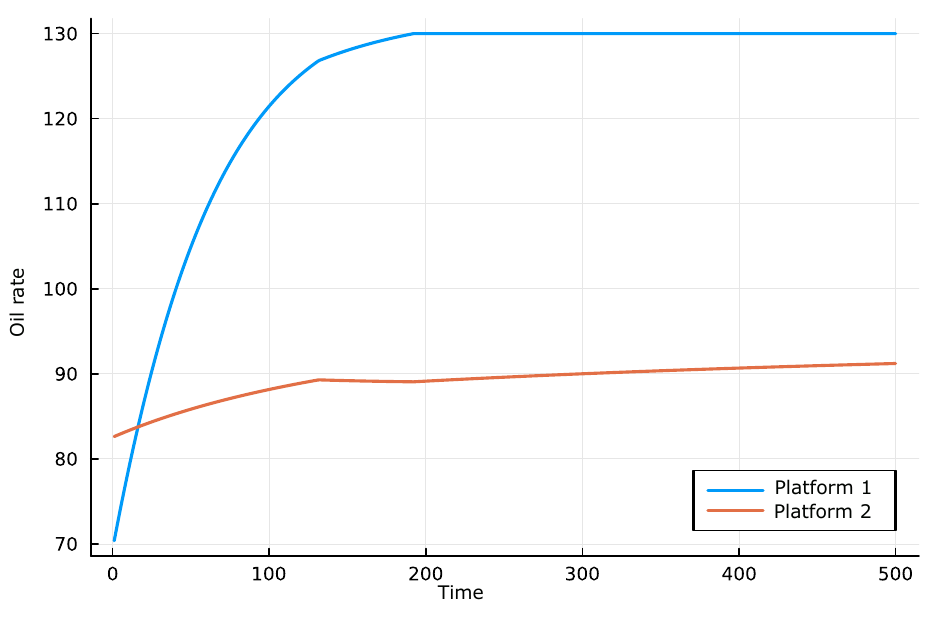}
         \caption{Outputs $y_i$ of the two platforms for the constrained case, in Sm$^3$day$^{-1}$}
         \label{fig:Outputs1_cstr}
     \end{subfigure}
     \hfill
          \begin{subfigure}[b]{0.3\textwidth}
         \centering
         \includegraphics[width=\textwidth]{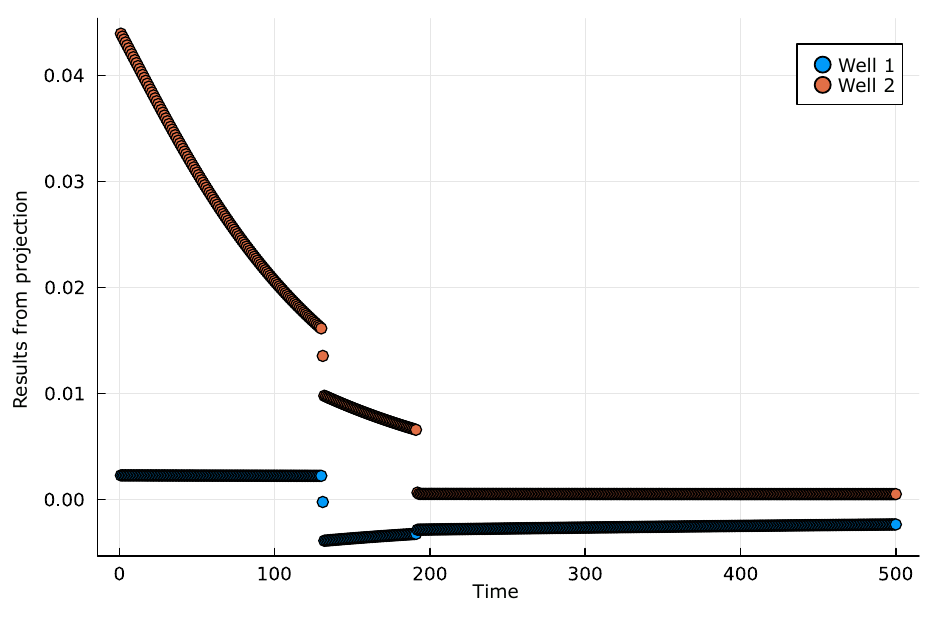}
         \caption{The solution to the projection over time}
         \label{fig:projection}
     \end{subfigure}
        \caption{Online Feedback Optimization for gas lift optimization over time (in iterations)}
        \label{fig:OFOGasLift}
\end{figure*}

\subsubsection{Sensitivity to gradient mismatch}
We use the formulas from Section \ref{sec:Sensitivity} to analyse the impact of the mismatch in the gradients $\nabla h$: $\nabla h^{\text{used}}(\mathbf{u}^k,\mathbf{y}^k)=\nabla h(\mathbf{u}^k,\mathbf{y}^k) + \beta_k$ where $\beta_k=\begin{bmatrix}\beta^{1}_k&\beta^{2}_k &0&0&0\\0&0&\beta_k^{3}&\beta_k^{4}&\beta_k^{5}\end{bmatrix}$. The additive mismatch was chosen to emulate  estimation error for every gradient separately. In particular, having a time-varying mismatch $\beta^i=[\beta_k^i]$ will allow us to analyse sensitivity if $\nabla h(\mathbf{u}^k,\mathbf{y}^k)\approx \overline{h}=$const, which is a practical approximation in OFO \cite{Tuning_Ortmann2024}.

\paragraph{Total sensitivity}

\begin{figure}[!tbp]
     \begin{subfigure}[b]{0.3\textwidth}
         \centering
         \includegraphics[width=\textwidth]{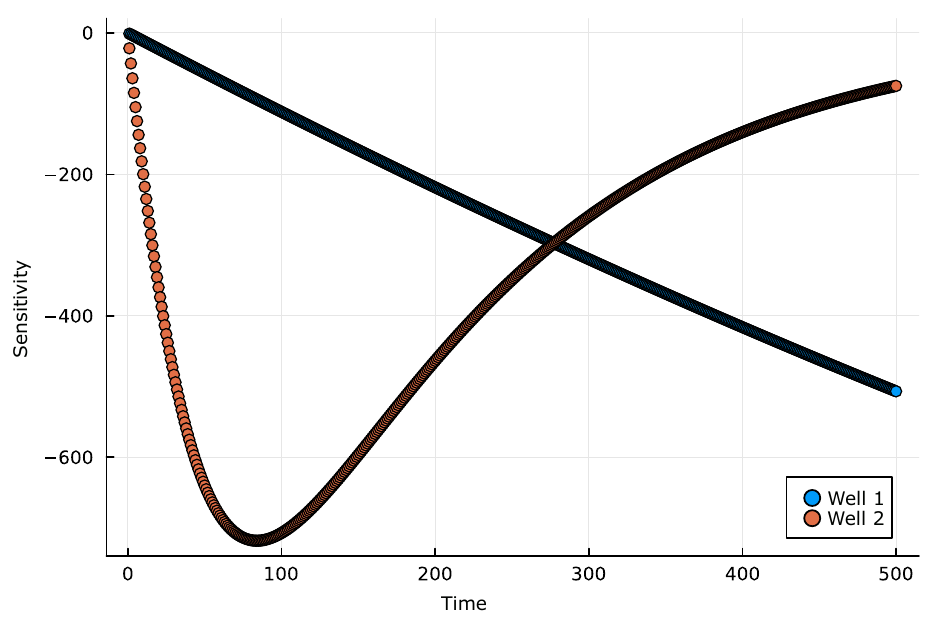}
         \caption{Total sensitivity of $\Phi$ at iteration $k$ to the mismatch $\beta^1$ and $\beta^2$ (Wells 1 and 2)}
         \label{fig:Derivatives1}
     \end{subfigure}
\hfill
     \centering
          \begin{subfigure}[b]{0.3\textwidth}
         \includegraphics[width=\textwidth]{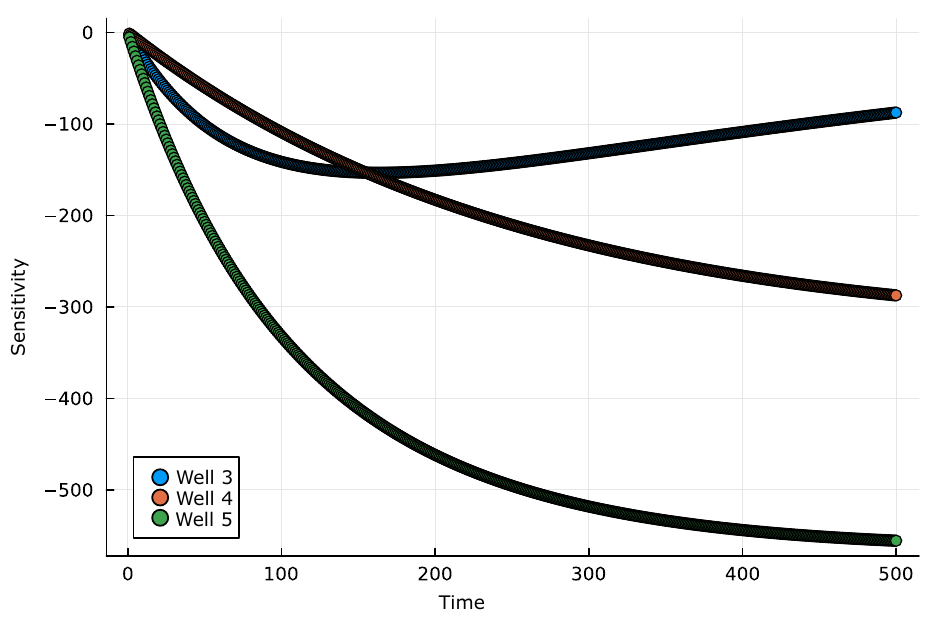}
        \caption{Total sensitivity of $\Phi$ at iteration $k$ to the mismatch $\beta^3$, $\beta^4$, $\beta^5$ (Wells 3, 4, and 5)}
        \label{fig:Derivatives2}
     \end{subfigure}
        \caption{Total sensitivity with respect to the mismatch in the five wells connected to two platforms, no coupling constraints}
        \label{fig:OFOSensitivity}
\end{figure}

The total sensitivities of the objective function at time $k$ are shown in Fig. \ref{fig:Derivatives1} for the first platform, and in Fig. \ref{fig:Derivatives1} for the second platform. The sensitivities to the mismatch are primarily influenced by how far from their respective optima the inputs are (Fig. \ref{fig:Inputs1} and Fig. \ref{fig:Characteristics}). The value of the objective is sensitive to $\beta^2$ in the initial period ($k<200$). This sensitivity is due to the ``steep'' characteristics of the second well (orange in Fig. \ref{fig:Characteristics}) around the initial guess ($u_2^0=7000$ Sm$^3$h$^{-1}$). This is because the platform with output $y_1$ reached the maximum of the well number two, bringing the derivative $\frac{\partial y_1}{\partial u_2}$ close to zero in \eqref{eq:PhiParamDeriv}, and thus reducing the impact of the mismatch on the objective. A similar sensitivity is observed for Well 3 (blue in Fig. \ref{fig:Derivatives2}) because Well 3 also reaches its optimum. However, as the characteristics of Well 3 are less ``steep'' than the characteristics of Well 1, the overall sensitivity is smaller. Conversely, Well 1, 4 and 5 are still operating away from their optima in Fig. \ref{fig:Characteristics}, so their respective sensitivities stay large. 

\begin{figure}[!tbp]
     \begin{subfigure}[b]{0.3\textwidth}
         \centering
         \includegraphics[width=\textwidth]{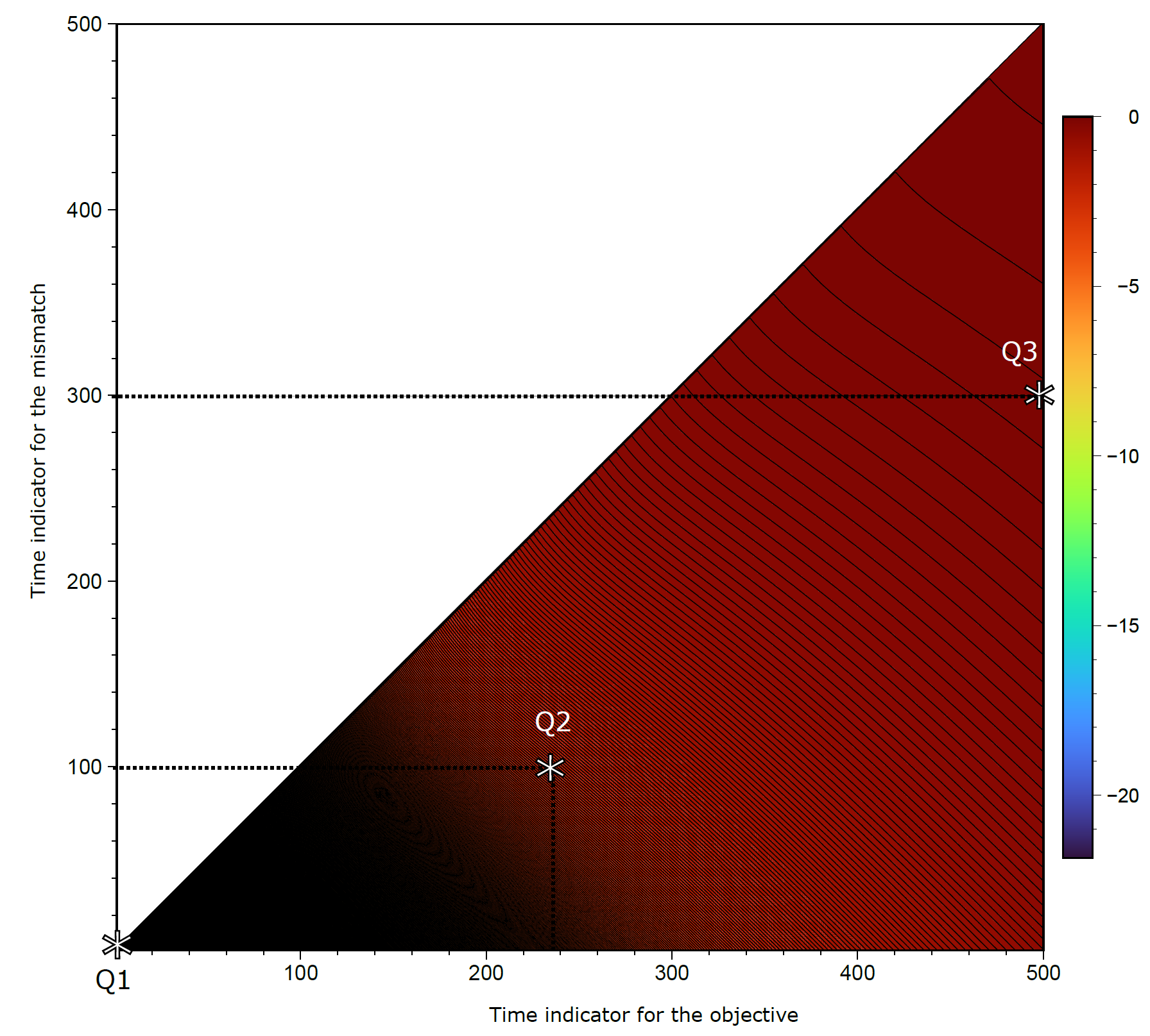}
         \caption{Instantaneous sensitivity of $\Phi$ to the mismatch $\beta$ in Well 2}
         \label{fig:Heatmap2}
     \end{subfigure}
\hfill
     \centering
          \begin{subfigure}[b]{0.3\textwidth}
         \includegraphics[width=\textwidth]{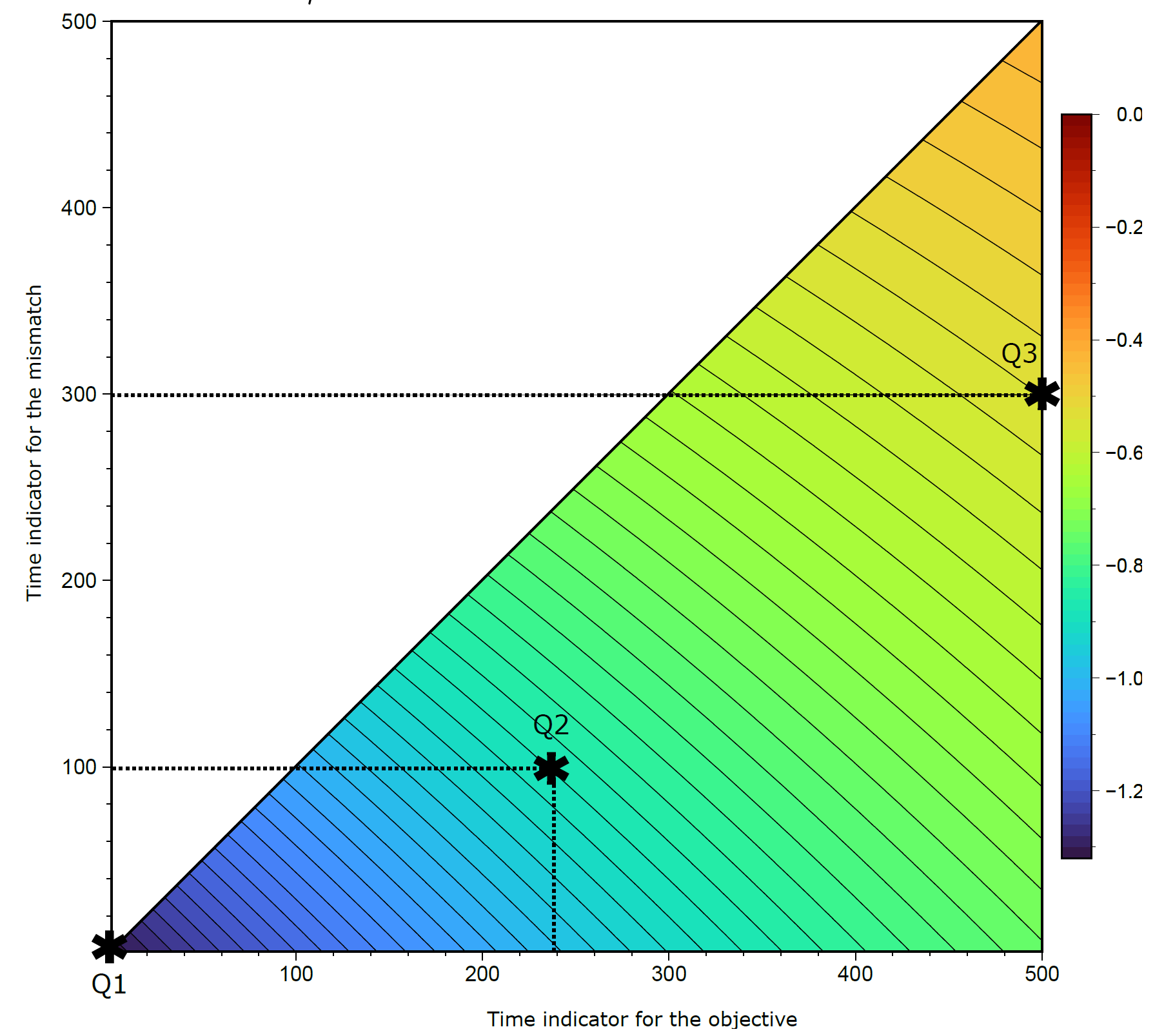}
        \caption{Instantaneous sensitivity of $\Phi$ to the mismatch $\beta$ in Well 4}
        \label{fig:Heatmap4}
     \end{subfigure}
        \caption{Instantaneous sensitivity with respect to the mismatch in Well 2 and 4, with three examples, Q1: $\partial\Phi^1/\partial\beta^i_1$, Q2: $\partial\Phi^{240}/\partial\beta^i_{100}$, Q3: $\partial\Phi^{500}/\partial\beta^i_{300}$, $i=2,4$, without coupling constraints}
        \label{fig:OFOHeatmap}
\end{figure}

\paragraph{Instantaneous sensitivity}
Figure \ref{fig:OFOHeatmap} show the instantaneous sensitivities of the objective with respect to $\beta^2$ and $\beta^4$, which were chosen as the most representative. In both cases, the largest absolute values are shown in dark blue (bottom left) and indicate that at the beginning of OFO, for small $k$, the objective changes under the influence of the mismatch. Conversely, for large $k$ (right hand-side of the figure), the impact of the mismatch at a single time step $j$ is small (yellow going into red). Such relationship between the objective function and the mismatch suggests that for large enough $k$ OFO is insensitive to mismatch in a single time step. The transition from dark blue to yellow and red for large $k$ also suggests that for large enough $k$ OFO becomes less sensitive to the initial mismatch. This in turn implies that iterative learning algorithms may be used for gradient estimation as long as their accuracy improves with iterations.

The results in Fig. \ref{fig:OFOHeatmap} allow us also to see the impact of the mismatch for individual wells. For instance, the mismatch affecting well 2 (Fig. \ref{fig:Heatmap2}) has a significantly larger impact at the beginning ($k$ close to zero, bottom left) than at the end of the chosen time horizon ($k$ close to 500, right). This result indicates that the characteristics of well 2 should be accurate to capture the impact of well 2 on the objective. Conversely, the instantaneous sensitivity of well 4 changes less from $k=0$ to $k=500$ (from left to right in Fig. \ref{fig:Heatmap4}). At the same time, we see that the sensitivity of the objective for $k=500$ depends on the timesteps, with a larger impact of mismatch in the initial time steps (light green, bottom right), getting close to zero with iterations (light orange, top right). Combining the transition from left to right and bottom to top suggests that while sensitivity to the mismatch in individual timesteps decreases, the final value is affected by mismatch in all past timesteps, in particular those far from the optimum.

\paragraph{Impact of coupling constraints}
The matrices in Section \ref{sec:NoCoupling} describe a case where the amount of gas available for injection to the platforms is unlimited, and the platforms have oversized capacity, which means that there is no coupling between the two wells.  As a special case, we now analyse the sensitivity to the mismatch if the amount of gas is limited, corresponding to an additional constraint $\sum\limits_{i=0}^5 u_i\leq 26000$ Sm$^{3}$day$^{-1}$, and the first platform has a decreased capacity, $y_1\leq 130$.

\begin{figure}[!tbp]
     \begin{subfigure}[b]{0.3\textwidth}
         \centering
         \includegraphics[width=\textwidth]{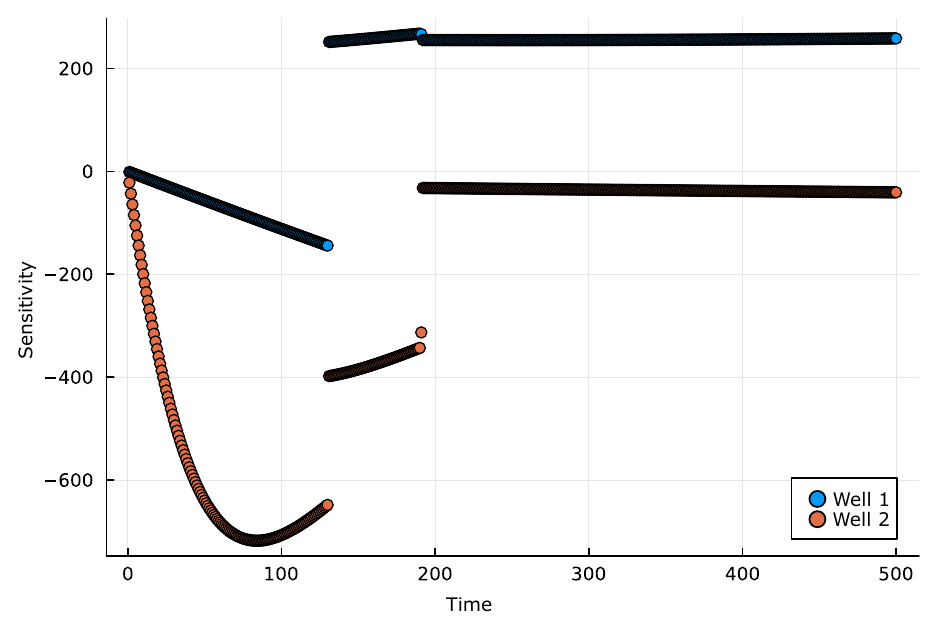}
         \caption{Total sensitivity of $\Phi$ at iteration $k$ to the mismatch $\beta^1$ and $\beta^2$ (Wells 1 and 2), for the case with coupling constraints}
         \label{fig:Derivatives1_cstr}
     \end{subfigure}
\hfill
     \centering
          \begin{subfigure}[b]{0.3\textwidth}
         \includegraphics[width=\textwidth]{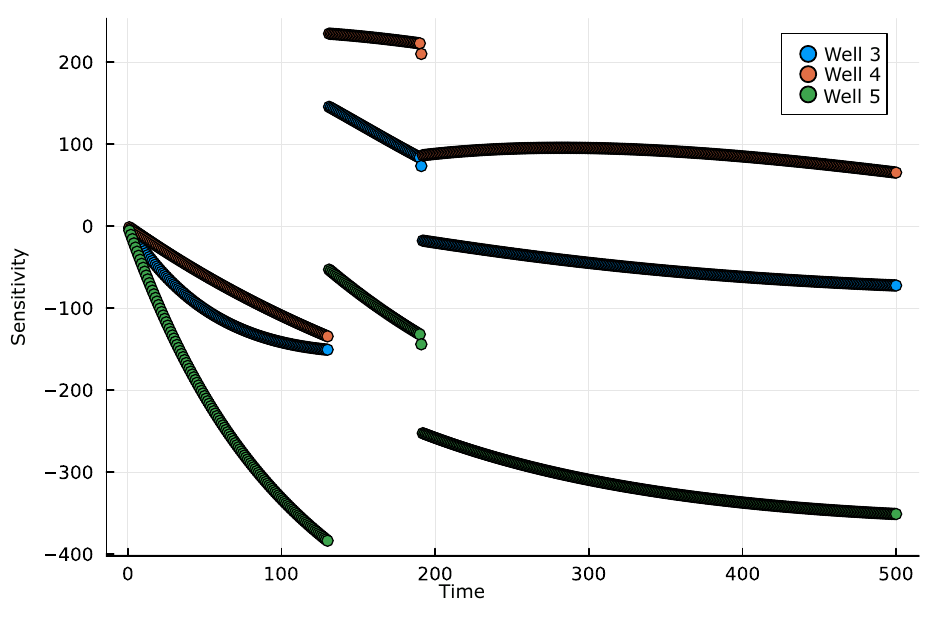}
        \caption{Total sensitivity of $\Phi$ at iteration $k$ to the mismatch $\beta^3$, $\beta^4$, $\beta^5$ (Wells 3, 4, and 5), for the case with coupling constraints}
        \label{fig:Derivatives2_cstr}
     \end{subfigure}
        \caption{Total sensitivity with respect to the mismatch in the five wells connected to two platforms, with coupling constraints}
        \label{fig:OFOSensitivity_cstr}
\end{figure}

The results of OFO for the case with coupling constraints are shown in the bottom row of Fig. \ref{fig:OFOGasLift} and the corresponding sensitivities are in Fig. \ref{fig:OFOSensitivity_cstr}. As long as no constraints are active, $k\leq 120$, the sensitivities in Fig.  \ref{fig:OFOSensitivity_cstr} are identical to those from Fig. \ref{fig:OFOSensitivity}. When the coupling constraint on the amount of gas becomes active, we note discontinuity in time is the sensitivity. This discontinuity is caused by the impact of active constraints on the solutions of the projection QP \eqref{eq:MismatchRewritten}, shown in Fig. \ref{fig:projection} for the first platform. The QP is solved at discrete iterations, so there is no assumption about continuity in time. Therefore, as the sensitivity of the objective depends on the sensitivity of all solutions in individual timesteps, we see the discontinuities in Fig. \ref{fig:OFOSensitivity_cstr}. A similar discontinuity in the sensitivities is observed at $k=200$ when the output constraint in the first platform becomes active, which is also caused by a discontinuity in the solutions from \ref{fig:projection}.

\paragraph{Impact of constant gradients}

The sensitivities $\frac{\partial\Phi}{\beta^i}$ shown in Fig. \ref{fig:Derivatives1} can be used to estimate how the objective will change if the mismatch changes. From the interpretation of total derivatives and treating $\Phi$ as a function of $\beta=[\beta^i]_{i=1,\ldots,5}$, we get:
\begin{equation}
\label{eq:approxDer}
    \Phi(\beta)-\Phi(0)\approx\sum\limits_{i=1}^5\frac{\partial\Phi}{\beta^i}\Delta\beta^i
\end{equation}
For the gas lift case study, we chose $\overline{h}=\nabla h(u_0,y(u_0))$ because the approximation at the initial condition was shown to work in practice \cite{Tuning_Ortmann2024}, and in consequence we got $\Delta\beta^i=[0,-0.04,-0.005,-0.001,-0.007]$, corresponding to the largest mismatch in absolute values. The maximal absolute value of the sensitivities is $[507,718,153,287,556]$. Plugging the values into \eqref{eq:approxDer} yields  $\Phi(\beta)-\Phi(0)\geq -32.6$ which means that we can expect a decrease of up to 14\% if the constant $\overline{h}$ is used. Because of taking the maximum values for the sensitivity and the mismatch, this is a conservative estimate, as also indicated in Fig. \ref{fig:Objective} where the actual decrease is 6.1\%. 

\section{Discussion}
\label{sec:Discussion}
The paper provides formulas for analysis of sensitivity of OFO. The results are validated in a synthetic case study, showing that the formulas can be used instead of finite differences to compute the sensitivities in both unconstrained and constrained optimization problems. The synthetic case study also indicated that the sensitivities of OFO depend on how long the OFO controller has been running. If the time was insufficient, the sensitivities of OFO increase as the optimum is not reached. Conversely, if OFO reaches the optimum, then it becomes less sensitive to changes in parameters. These results suggest that future work could focus on  analysis of sensitivity with respect to time.

The second numerical case study presents an application of sensitivity analysis with respect to the model mismatch in a gas lift optimization problem. The case study further demonstrated that the sensitivity of OFO depends on how long the controller was running. In particular, the sensitivities to the model mismatch in a single time step decrease with time, which suggests that iterative algorithms with accuracy improving with iterations may be used for gradient estimation while preserving the optimum.

The results for the case with coupling constraints emphasize the dependence on time in OFO. The formulas in Section \ref{sec:Sensitivity} were derived by setting a constant final time $T_F$ and treating OFO as a system of equations parametrised by $\mathbf{p}$:
\begin{equation} \label{eqn:Equations}
\begin{aligned}
\Phi^{T_F}= &{}\Phi(u^{T_F},y^{T_F},\mathbf{p})\\
    y^0 = &{}h(u^0)\\
    u^{1} = &{}u^0 + \alpha\widehat{\sigma}_\alpha (u^0,y^0,\mathbf{p}) \quad \\
    \vdots\\
    y^{k} =&{} h(u^{k})\\
    u^{k+1} = &{}u^k + \alpha\widehat{\sigma}_\alpha (u^k,y^k,\mathbf{p})\\
    \vdots\\
    y^{T_F-1} = &{}h(u^{T_F-1})\\
    u^{T_F} =&{} u^{T_F-1} + \alpha\widehat{\sigma}_\alpha (u^{T_F-1},y^{T_F-1},\mathbf{p})
\end{aligned}
\end{equation}
where $\widehat{\sigma}_\alpha (u^k,y^k,\mathbf{p})$ is the solution to \eqref{eq:MismatchRewritten} at iteration $k$. This assumption is equivalent to discretizing a continuous gradient flow algorithm with a constant time step \cite{Non_Haeberle2020}. However, \eqref{eqn:Equations} can be equivalently summarised as a nonlinear equation:
\begin{equation}
    F(\Phi^{T_F},\mathbf{u},\mathbf{y},\mathbf{p})=0
\end{equation}
treating all arguments as independent variables, with no explicit dependence on time. A natural extension to the sensitivity analysis will now be to consider the impact of time, both as time step and the final time.

\section{Conclusions and future works}
\label{sect:concl}

The objective of the paper is to facilitate the analysis of Online Feedback Optimization controllers with respect to their parameters. The importance of analysing the impact of the parameters of an optimization-based controller has been shown and used for traditional controllers, such as Model Predictive Control, but the analysis for Online Feedback Optimization remains under-explored. This paper addresses this gap by providing closed-form expressions for the sensitivity of Online Feedback Optimization to its parameters. 

In the future, the sensitivities can be used for tuning of the parameters of the controller, as well as checking robustness to model mismatch. Future work will also include the analysis of sensitivity with respect to time.

%

\end{document}